\documentclass[12pt]{article}
\usepackage{tikz}
\usetikzlibrary{arrows.meta}
\usepackage{amssymb}

\begin{document}
\topmargin 0.0in
\textheight 8.5in
\textwidth 5in
\oddsidemargin 0.5in
\title{  Cardy's formula does not hold on some 2D  lattices for  critical two-dimensional percolation
} 
\author{Yu Zhang}
\baselineskip .20in
\maketitle
\begin{abstract}
The scaling limit of crossing probabilities  is believed to satisfy a conformal mapping  formula, called Cardy's formula,  in two-dimensional percolation at the criticality.  The formula has been confirmed to hold for site percolation on the equilateral triangular lattice. In this paper, we show that   Cardy's formula could not hold for some two-dimensional  triangular and square-type lattices, in particular for some periodic 2D graphs.

\end{abstract}
\section{ Introduction and statement of results.}
Percolation is one of the fundamental stochastic models studied by mathematicians. Since percolation is
one of the simplest models that exhibit a phase transition, it also becomes  
one of statistical physicists' favorite models for studying critical phenomena. 
In particular, the case of two dimensions is very
special since it links the conformal field theory with the singularity of a phase
transition at the criticality.  
We first introduce a two-dimensional graph ${\cal G}=\{V, E\}$ for vertex and edge sets $V$ and $E$ with the following three requirements (see pages 10-12 in Kesten (1982)):\\
(1) There exists  a $k< \infty$, such that there are at most $k$ edges of ${\cal G}$ incident to any vertex of ${\cal G}$.\\
(2) All edges of ${\cal G}$ have a finite diameter. Every compact set of ${\bf R}^2$ intersects only finitely many edges of ${\cal G}$.\\
(3) ${\cal G}$ is connected where two vertices $u$ and $v$ are connected when there is an edge with two vertices $u$ and $v$. A sub-graph is  connected if all of its two vertices  are connected.\\
We say ${\cal G}$ is a {\em periodic graph} if (1)--(3) hold and
${\cal G}$ is imbedded in ${\bf R}^2$ in such a way that each coordinate vector of ${\bf R}^2$ is a period for the image.
A period for image means that  $v\in {\bf R}^2$ is a vertex of ${\cal G}$ if and only if $v+ k_1\xi_1+k_2 \xi_2$ is a vertex of ${\cal G}$ for all integers $k_i$ for $i=1,2$ in ${\bf Z}$, where $\xi_i$  denotes the $i$-th coordinate vector of ${\bf R}^2$. 

 One of the most popular 2D lattices is  the square lattice ${\cal S}_\delta=\{{\bf Z}^2_\delta, {\cal E}_\delta ^2\}$  with mesh $\delta>0$.
 More precisely, the square lattice  has   vertex set  
$${\bf Z}^2_\delta=\{(\delta i, \delta j): \mbox{ $i$ and $j$ are integers and $\delta>0$}\},$$  
and  edge set  ${\cal E}_\delta^2$ connecting a
pair of vertices
$u=(u_1, u_2)$ and $v=(v_1, v_2)$ 
with $d(u,v)=\delta$, where $d(u,v)$ is the Euclidean distance between $u$ and $v$. 
 Graph ${\cal S}_\delta$  is a periodic graph.
For any two sets ${A}$ and ${B}$, we also define their distance by
$$d({A}, {B})= \min\{ d(u, v): u\in {A}, v\in {B}\}.$$
Note that vertices connected by  edges consist of squares. We call these squares
$\delta$-{\em squares}.
When $\delta=1$, ${\cal S}_1$ is the  standard square lattice originally considered  by Broadbent and Hammersley (1957).

We also  introduce the triangular lattice with mesh $\delta$: Divide ${\bf R}^2$ into isosceles triangles by means of the horizontal lines
and lines under an angle of $0< \theta< \pi/2$ and $\pi-\theta$, with the first coordinate axis through the points $(\delta k, 0), k\in {\bf Z}$ 
(see Figure 1). The vertices of the graph are the vertices of the  isosceles triangles, and the edges are the line segments connecting two vertices
of the same triangle. We divide these line segments into {\em horizontal, north-east, } and {\em north-west} edges
with vertices $u=(u_1, u_2)$ and $v=(v_1, v_2)$  such that $d(u,v)=\delta$ for horizontal edges, and $d(u, v)=k\delta$ for hypotenuse edges.
  In particular, if  $k=1$ and $\delta =1$, then it is the standard equilateral triangular lattice.
We denote the above graph by ${\cal T}_\delta(k)=\{{\bf V}_\delta(k), {\bf T}_\delta(k)\}$ for vertex set ${\bf V}_\delta(k)$ and edge set ${\bf T}_\delta(k)$, respectively. 
 The standard equilateral triangular lattice ${\cal T}_\delta (1)={\cal T}_\delta =\{ {\bf  V}_\delta, {\bf T}_\delta\}$ is not a periodic graph. In fact,  $(1,0)$ is a period for this imbedding, but $(0,1)$ is not.
 However, ${\cal T}_\delta({1/\sqrt{3}})$ is a periodic graph.   

 We consider the square lattice ${\cal  S}_\delta$ and add  a north-east edge in each $\delta$-square. We denote the graph by 
 ${\cal S}_{\delta}^{ne}= \{{\bf Z}_\delta^2, {\cal E}_\delta ^2 (ne)\}$ for the vertex and edge sets (see Figure. 2). In fact, if we rotate ${\cal S}_{\delta/\sqrt{2}}^{ne}$ into
 $\pi/4$ counterclockwise, it is the triangular lattice ${\cal T}_\delta(1/\sqrt{2})$. ${\cal  S}_\delta^{ne}$ is a periodic lattice.

Now we introduce percolation model on graph ${\cal G}$ (see Kesten (1982)).
Each vertex in ${\cal G}$ is either {\em open} or {\em closed}
 independently with  probability $p$ and $1-p$.
This model is called {\em site percolation} on ${\cal G}$.
The corresponding probability measure on the configurations of open 
and closed  vertices  is denoted by ${\bf P}_{p, \delta}$. 
A path  from $u$ to $v$ for $u, v\in {\cal G}$ is a sequence $(v_0,... ,v_{i}, v_{i+1},... ,v_n)$
with distinct vertices $v_i$, $0\leq i\leq n$,
 and $v_0=u$ and $v_n=v$  such that $v_i$ and $v_{i+1}$
are connected. 
If all the vertices in a path are open or closed, then the path is called an open or a closed path, respectively.
The open cluster of the vertex $x$, 
${\bf C}(x)$,
consists of all vertices that are connected to $x$ by an open path.
For any collection ${A}$ of vertices, $| A |$  
denotes the cardinality of $A$. We choose ${\bf 0}$ as the origin. The percolation probability is
\begin{eqnarray*}
\theta (p)= {\bf P}_{p,\delta}(|{\bf C}({\bf 0})|=\infty),
\end{eqnarray*}
and the critical probability is 
$$p_c=\sup\{p:\theta (p)=0\}.$$
It is well known (see Kesten (1982)) that
$ 1> p_c > 0.5$ for site percolation on ${\cal S}_\delta$ and that $ p_c=0.5$ for site percolation on ${\cal T}_{\delta}$.

\section{ Conformal mapping for crossing probabilities and Cardy's formula.}
Now we will consider the complex plane ${\bf C}$.
Let ${D} \subset {\bf C}$ be a simple connected, bounded domain, whose boundary $\partial {D}$ is a Jordan curve. ${D}$ is called 
a {\em Jordan domain}. We select four points $a, x, b, c$ from $\partial D$ with $x\in \overline{ab}$ such that 
$$\partial D=\overline{ax}\cup\overline{xb}\cup\overline{bc}\cup\overline{ca}.$$
We denote the {\em crossing probability} by 
$$\pi_{p,\delta}({D}, \overline{ax}, \overline{bc})= {\bf P}_{ p,\delta}\left(\exists \mbox{ an open $\bullet$-path in $D$ from $\overline{ax}$ to $\overline{bc}$}\right).\eqno{}$$
For simplicity,  if $p=p_c$, let
$$\pi_{\delta}({D}, \overline{ax}, \overline{bc})=\pi_{p_c,\delta}({D}, \overline{ax},\overline{bc})={\bf P}_{\delta}\left(\exists \mbox{ an open $\bullet$-path in $D$ from $\overline{ax}$ to $\overline{bc}$}\right).\eqno{}$$
With these definitions, it is known (see Kesten (1982)) that for any Jordan domain ${D}$ and $a, x, b, c$,
$$\lim_{\delta\rightarrow 0} \pi_{p,\delta}({D}, \overline{ax}, \overline{bc})=\left \{ \begin{array}{ll}0 & \mbox{ if $p < p_c,$}\\
1 &\mbox{ if $p>p_c.$}
\end{array}\right. \eqno{(2.1)}
$$
When $p=p_c$, it is also known   (see Kesten (1982)) that
$$0<\liminf_{\delta \rightarrow 0} \pi_{\delta}({D}, \overline{ax},\overline{bc})\leq \limsup_{\delta \rightarrow 0} 
 \pi_{\delta }({D}, \overline{ax},\overline{bc})<1.\eqno{(2.2)}$$
With (2.2),  Langlands et al. (1994)  conjectured  
the existence of the scaling limit:
$$\lim_{\delta \rightarrow 0} \pi_{\delta}({D}, \overline{ax},\overline{bc})
=\pi({D}, \overline{ax},\overline{bc}).\eqno{(2.3)}$$
  By the Riemann mapping theorem, there is a unique conformal mapping $\phi$ that maps ${D}$ to the unit equilateral triangle $\Delta ABC$ on ${\bf C}$ with 
$$\phi(a)=A=0, \phi(b)=B=1,   \phi(c)=C=e^{i\pi/3}, \phi(x)=X,\eqno{(2.4)}$$
for $X\in \overline{AB}$.   

With this map, it was  believed  by Cardy (see Cardy (1992) and Grimmett (1999 p. 346))
that the limit probability of (2.3) exists and is equal to $X$ for most graphs.
This conjecture is called Cardy's formula with Carleson's version.  \\

{\bf Cardy's formula.} {\em If $D$ is a Jordan domain with four points $a,x,b,c$  for any $x\in \overline{ab}$ defined  above, and $\phi$ is the Riemann mapping such that $ \phi(a)=A=0, \phi(b)=B=1,   \phi(c)=C=e^{i\pi/3}, \phi(x)=X$, then }
$$\lim_{\delta \rightarrow 0} \pi_\delta (D, \overline {ax}, \overline {bc})=X.$$

Smirnov's  celebrated  work (2001) showed that the scaling limit in (2.3) indeed exists
and satisfies  Cardy's formula on ${\cal T}_{\delta}$.   With his work, the power laws and  the existence of critical exponents with the scaling relations, 
one of the most important questions in physics, are proved to exist in the percolation model.
The following theorem shows, however, that Cardy's formula may not  hold for some graphs.\\ 

{\bf Theorem.} {\em If $k\neq 1$, then Cardy's formula does not hold for site percolation on ${\cal T}_{\delta}(k)$.
Also, Cardy's formula does not hold on ${\cal S}_\delta^{ne}$.}\\

${\cal T}_\delta$ has only three kinds of edges: north-east, north-west, and horizontal edges. We denote  by 
${\bf T}_\delta^{ne}$, ${\bf T}_\delta^{nw}$, and ${\bf T}_\delta^h$ the edge sets with north-east and horizontal edges,   with north-west and horizontal edges,  and with north-east and north-west edges, respectively.  Let
${\cal T}_\delta^{ne}=\{{\bf V}_\delta, {\bf T}_\delta^{ne}\}$, ${\cal T}_\delta^{nw}=\{{\bf V}_\delta, {\bf T}_\delta^{nw}\}$, and ${\cal T}_\delta^{h}=\{{\bf V}_\delta, {\bf T}_\delta^{h}\}$ be the three corresponding graphs.  By using the method of the theorem, if one can show that Cardy's formula  holds on one of the above three graphs, then Cardy's formula will not hold on the square lattice ${\cal S}_\delta$.  On the other hand, if one can show that Cardy's formula holds on the square lattice,  then Cardy's formula will not hold on the above three graphs.
We believe the first statement, since Cardy's formula holds on ${\cal T}_\delta$ but not on square-type lattice ${\cal S}_\delta^{ne}$ in the theorem.\\

{\bf Conjecture.} {\em Cardy's formula only holds on ${\cal T}_\delta$, ${\cal T}_\delta^{ne}$, ${\cal T}_\delta^{nw}$, and
${\cal T}_\delta^{h}$. }\\

{\bf Remark.} We believe that 
the scaling limit exists for most graphs at least for all periodic graphs even though the limit may not satisfy Cardy's formula. 
In fact, by Smirnov's result (2001), one can show that the scaling limit in (2.3) exists  for the graphs  ${\cal T}_{\delta}(k)$ for all $k>0$, but Cardy's formula does not hold.  If the limit exists, by a standard computation (see Werner (2008)),  the exploration process  for ${\cal T}(k)$ converges.  One may ask whether or not  the exploration process  for ${\cal T}(k)$ converges to $SLE_\kappa$ for all $k$? If it converges, then what is $\kappa=\kappa(k)$? 
If one can show Cardy's formula holds on one  of graphs ${\cal T}_\delta^{ne}$, ${\cal T}_\delta^{nw}$, or ${\cal T}_\delta^{h}$,
 then by a standard computation (see Werner (2008)),
the power laws, the existence of critical exponents with the scaling relations, and universality are proved to hold on the square lattice ${\cal S}_\delta$, even though Cardy's formula does not hold on ${\cal S}_\delta$.

\section{Proof of Theorem.}
\begin{figure}
\begin{center}
\setlength{\unitlength}{0.0125in}%
\begin{picture}(200,220)(67,680)
\thicklines
\put(-50,680){\line(1,3){100}}
\put(-10,680){\line(1,3){80}}
\put(30,680){\line(1,3){60}}
\put(70,680){\line(1,3){40}}
\put(110,680){\line(1,3){20}}

\put(150,680){\line(-1,3){100}}
\put(110,680){\line(-1,3){80}}
\put(70,680){\line(-1,3){60}}
\put(30,680){\line(-1,3){40}}
\put(-10,680){\line(-1,3){20}}
\put(-50,680){\line(1,0){200}}
\put(-30,740){\line(1,0){160}}
\put(-10,800){\line(1,0){120}}
\put(10,860){\line(1,0){80}}
\put(30,920){\line(1,0){40}}

\put(170,680){\line(2,3){100}}
\put(210,680){\line(2,3){80}}
\put(250,680){\line(2,3){60}}
\put(290,680){\line(2,3){40}}
\put(330,680){\line(2,3){20}}

\put(370,680){\line(-2,3){100}}
\put(330,680){\line(-2,3){80}}
\put(290,680){\line(-2,3){60}}
\put(250,680){\line(-2,3){40}}
\put(210,680){\line(-2,3){20}}

\put(170,680){\line(1,0){200}}
\put(190,710){\line(1,0){160}}
\put(210,740){\line(1,0){120}}
\put(230,770){\line(1,0){80}}
\put(250,800){\line(1,0){40}}

\put(70,920){\circle*{8}}
\put(50,860){\circle*{8}}
\put(30,800){\circle*{8}}
\put(50,740){\circle*{8}}
\put(30,680){\circle*{8}}

\put(290,800){\circle*{8}}
\put(270,770){\circle*{8}}
\put(250,740){\circle*{8}}
\put(270,710){\circle*{8}}
\put(250,680){\circle*{8}}

\put(-55,660){\mbox{$\alpha=0$}{}}
\put(125,660){\mbox{$\beta=1$}{}}
\put(70,670){\mbox{$x$}{}}
\put(290,670){\mbox{$x$}{}}
\put(-45,680){\mbox{$\kappa$}{}}

\put(170,660){\mbox{$A=0$}{}}
\put(360,660){\mbox{$B=1$}{}}
\put(60,975){\mbox{$\gamma=2e^{i\pi/3}$}{}}
\put(260,835){\mbox{$C=e^{i\pi/3}$}{}}
\end{picture}
\end{center}
\caption{ {\em The left graph is $\bigtriangleup \alpha\beta\gamma$ with sites on  ${\cal T}_{\delta}( 2 )$, and the right is $\bigtriangleup  ABC$ with sites on ${\cal T}_{\delta}$.  Each open path, shown as a solid circle path,  from $\overline{\beta \gamma}$ to $\overline{\alpha x}$ corresponds  to an open path from $\overline{BC}$ to $\overline{Ax}$ with the same probability.}}
\end{figure}
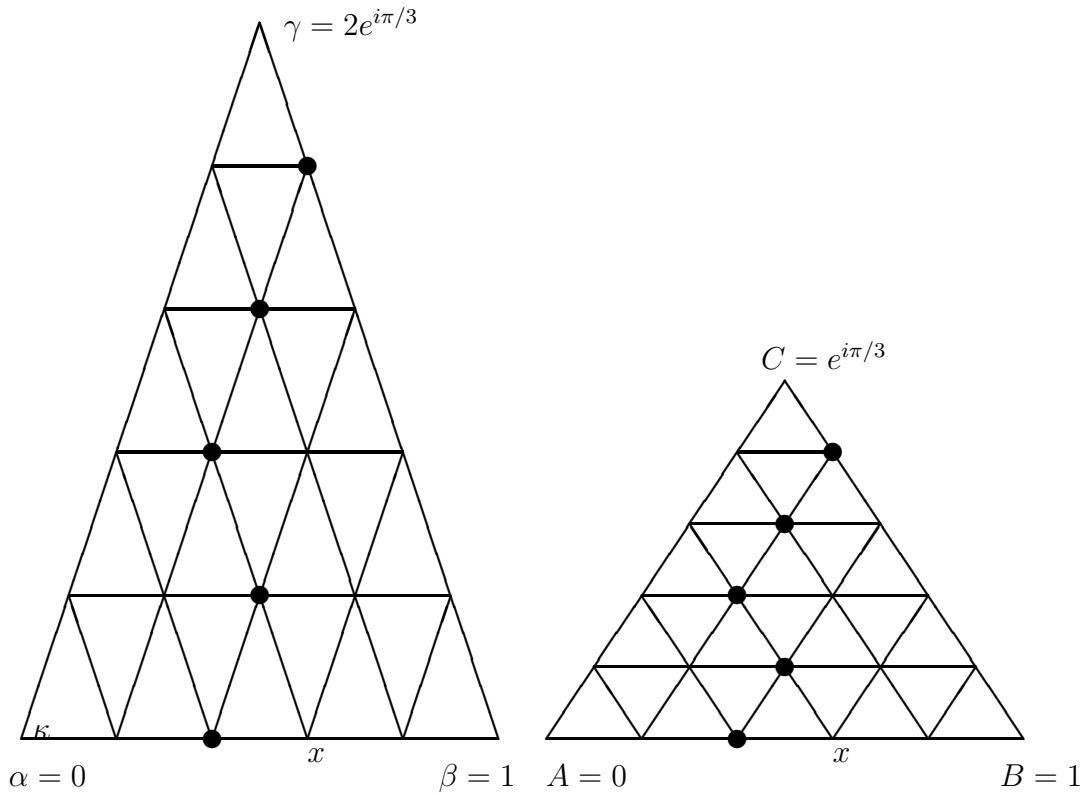

 Without loss of generality, we consider ${\cal T}_{\delta}(2)$ and triangle $ \bigtriangleup\alpha\beta\gamma$ with
$\alpha=0, \beta=1$, and $\gamma=2e^{i\pi /3}$ (see Figure 1). 
In addition, we also consider ${\cal T}_{\delta}$ and triangle $ \bigtriangleup ABC$ with
$A=0, B=1$, and $C=e^{i\pi /3}$. 
We select  a vertex $x\in {\cal T}_{\delta}(2)\cap \overline{\alpha\beta}$. We also select the same vertex $x\in {\cal T}_{\delta}\cap \overline{AB}$.
Note that each triangle of ${\cal T}_{\delta}$  is stretched vertically two times to be a triangle of  ${\cal T}_{\delta}(2)$, so  any corresponding two open paths from $x$ to $\overline{\beta\gamma}$ in ${\cal T}_{\delta}(2)$ and from $x$ to $\overline{BC}$ in ${\cal T}_{\delta}$ have the same probability (see Figure 1). In other words,
$$ \pi_{\delta}({\bigtriangleup\alpha\beta\gamma}, \overline{\alpha x},\overline{\beta\gamma})= \pi_{\delta}({\bigtriangleup ABC}, \overline{A x},\overline{BC}).\eqno{(3.1)}$$
Thus, by Smirnov's theorem,
$$\lim_{\delta\rightarrow 0} \pi_{\delta}({\Delta\alpha\beta\gamma}, \overline{\alpha x},\overline{\beta\gamma})
=\pi({\bigtriangleup \alpha\beta\gamma}, \overline{\alpha x},\overline{\beta\gamma})=\pi({\bigtriangleup ABC}, \overline{A x},\overline{BC})=x.\eqno{(3.2)}$$
Note that the unique conformal map from $\Delta \alpha\beta\gamma$ to $\Delta ABC$ cannot be linear.  We can show this by using the following precise  maps.
Let $\angle \gamma\alpha \beta=\kappa$.
Schwarz-Christoffel's  transformation $\phi_1(w)=z $ maps  from the upper half-plane to $\Delta \alpha\beta\gamma$ such that $\phi_1(0)=\alpha$, $\phi_1(1)=\beta$, and $\phi_1(\infty)=\gamma$ with
$$x=\phi_1(w) =A_1\int^{w} _0 y^{\kappa/\pi-1}(1-y)^{\kappa/\pi-1}dy\mbox{ for }w\in(0, 1) \eqno{(3.3)}$$
for a constant $A_1$. In addition, 
Schwarz-Christoffel's  transformation $\phi_2(w)=z $ maps from the upper half-plane to $\Delta ABC$ such that
$\phi_2(0)=A$, $\phi_1(1)=B$ and $\phi_1(\infty)=C$ with
$$X=\phi_2(w) =A_2\int^{w} _0 y^{-2/3}(1-y)^{-2/3}dy\mbox{ for }w\in(0, 1)\eqno{(3.4)}$$
 for a constant $A_2$. By (3.2) and (3.3),
$${dx\over dw}=A_1w^{\kappa/\pi-1}(1-w)^{\kappa/\pi-1}\mbox{ for }w\in (0, 1),\eqno{(3.5)}$$
and 
$${dX\over dw}=A_2 w^{-2/3}(1-w)^{-2/3}\mbox{ for }w\in (0, 1).\eqno{(3.6)}$$
If Cardy's formula holds, then by (3.2),
$$X=x.\eqno{(3.7)}$$
By (3.5)--(3.7),
$$A_1w^{\kappa/\pi-1/3}(1-w)^{\kappa/\pi-1/3}=A_2 .\eqno{(3.8)}$$
Since $\kappa > \pi/3$,  (3.8) cannot hold if $w\downarrow 0$ or $w\uparrow 1$.
So Cardy's formula does not hold for ${\cal T}_{\delta}( 2)$.   

We may consider $ \bigtriangleup \alpha\beta\gamma$ for $\alpha=0$, $\beta=1$, and $\gamma=k e^{i\pi/3}$ for any positive $k\neq 1$. 
 Let $\phi$ be the conformal map that maps $\bigtriangleup \alpha\beta\gamma$ to $\bigtriangleup ABC$. By the same proofs in (3.1)--(3.8),  if
 $x\in (0, 1)$, then
 $$\phi(x)\neq x.\eqno{(3.9)}$$
  
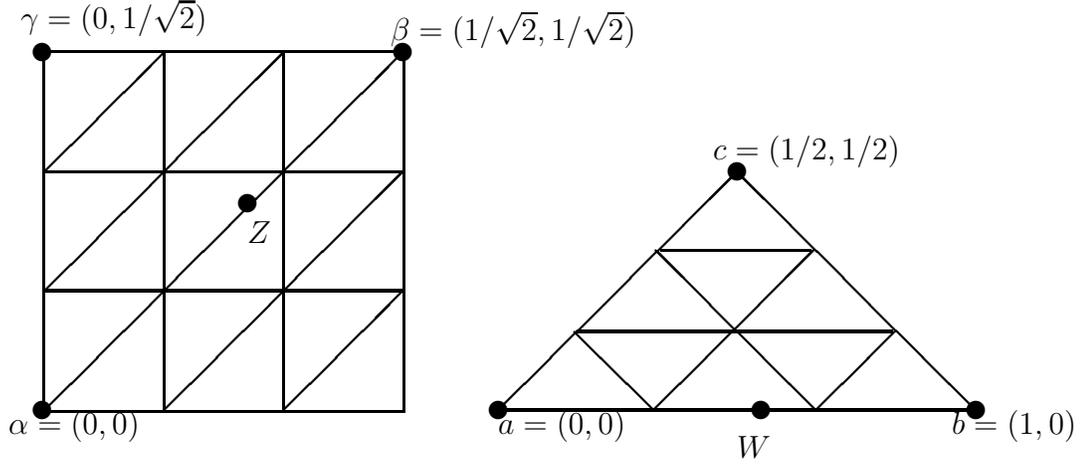
\begin{figure}
\begin{center}
\setlength{\unitlength}{0.0125in}%
\begin{picture}(200,100)(67,680)
\thicklines

\put(-20, 680){\framebox(150,150){}}
\put(-20, 730){\line(1,0){150}}
\put(-20, 780){\line(1,0){150}}
\put(30, 680){\line(0,1){150}}
\put(80, 680){\line(0,1){150}}

\put(-20,680){\line(1,1){150}}
\put(-20,730){\line(1,1){100}}
\put(-20,780){\line(1,1){50}}
\put(30,680){\line(1,1){100}}
\put(80,680){\line(1,1){50}}

\put(170,680){\line(1,1){100}}
\put(235,680){\line(1,1){67}}
\put(303,680){\line(1,1){33}}

\put(370,680){\line(-1,1){100}}
\put(303,680){\line(-1,1){67}}
\put(235,680){\line(-1,1){33}}

\put(170,680){\line(1,0){200}}
\put(203,713){\line(1,0){132}}
\put(238,747){\line(1,0){63}}

\put(170,680){\circle*{8}}
\put(270,780){\circle*{8}}
\put(370,680){\circle*{8}}
\put(-21,830){\circle*{8}}
\put(-21,680){\circle*{8}}
\put(130,830){\circle*{8}}

\put(125,835){\mbox{$\beta=(1/\sqrt{2},1/\sqrt{2})$}{}}
\put(-30,840){\mbox{$\gamma=(0,1/\sqrt{2})$}{}}
\put(65,750){\mbox{$Z$}{}}
\put(270,660){\mbox{$W$}{}}
\put(65,767){\circle*{8}}
\put(280,680){\circle*{8}}
\put(-35,670){\mbox{$\alpha=(0,0)$}{}}

\put(170,670){\mbox{$a=(0,0)$}{}}
\put(360,670){\mbox{$b=(1,0)$}{}}
\put(260,785){\mbox{$c=(1/2, 1/2)$}{}}
\end{picture}
\end{center}
\caption{ {\em The left graph is ${\cal S}_{\delta/\sqrt{2}}^{ne}$. After rotation $\pi/4$ counterclockwise, $\bigtriangleup \alpha\beta\gamma$ is $\bigtriangleup abc$ in graph ${\cal T}_{\delta} (1/\sqrt{2})$ (right graph).  $Z$ on $\overline{\alpha \beta}$ will be $W$ on $\overline{ab}$.}}
\end{figure}

Now we show that Cardy's formula does not hold on ${\cal S}^{ne}_\delta$.  In fact, we only need to show that Cardy's formula does not hold on ${\cal S}^{ne}_{\delta/\sqrt{2}}$. To begin, we assume that Cardy's formula holds on ${\cal S}^{ne}_{\delta/\sqrt{2}}$.
Let $\bigtriangleup \alpha \beta\gamma$ be the triangle with vertices $\alpha=(0,0)$, $\beta=(1/\sqrt{2}, 1/\sqrt{2})$, and $\gamma=(0,1/\sqrt{2})$ (see Figure 2).  We divide $\overline{\alpha\beta}$ into $\delta^{-1}/\sqrt{2}$ many intervals. For each $Z\in \bigtriangleup \alpha \beta\gamma$, we use the rotation map $W= \phi_1(Z)=e^{-i\pi/4} Z$ to map 
$\bigtriangleup \alpha \beta\gamma$ into $\bigtriangleup abc$ for $a=0$, $b=1$, and $c=1/2+i/2$.  Thus, 
$$\phi_1({\cal S}^{ne}_{\delta/\sqrt{2}})= {\cal T}_{\delta}( 1/\sqrt{2}).\eqno{(3.10)}$$
We select $Z \in \overline{\alpha\beta}$ (see Figure 2). Note that  $\bigtriangleup \alpha \beta\gamma$ can be viewed to be  $\bigtriangleup abc$
after rotating $\bigtriangleup abc$ at the origin $\pi/4$ counterclockwise. 
Note also that  (3.1) works for any $k$,  so if Cardy's formula holds, then
$$\pi(\bigtriangleup \alpha\beta\gamma, \overline{\alpha Z}, \overline{\beta\gamma})=d(0, Z).\eqno{(3.11)}$$
Since $\phi_1$ is a  rotation map,  by (3.11),
$$\phi_1(Z)=W=d(0, Z).\eqno{(3.12)}$$
Let $\phi_2$ be  Schwarz-Christoffel's  transformation that maps from $\bigtriangleup abc$ to the unit equilateral triangle $\Delta ABC$.
By (3.9), this Schwarz-Christoffel's  transformation  is not an identity map, so
$$\phi_2(W) \neq W.\eqno{(3.13)}$$
On the other side, note that $\phi_2\phi_1$ is the conformal mapping from $\bigtriangleup \alpha \beta\gamma$ to $\bigtriangleup ABC$, so if Cardy's formula holds, then  
$$\phi_2(W)=\phi_2\phi_1(Z)=d(0, Z)=W.\eqno{(3.14)}$$
Therefore, (3.13) and (3.14) cannot hold together, so Cardy's formula does not hold on ${\cal S}^{ne}_\delta$.  
The theorem follows. $\blacksquare$

\begin{center}{\large \bf References} \end{center}
Broadbent, S. R. and  Hammersley, J. M. (1957). Percolation processes I. Crystals and mazes. {\em Math. Proc. of the Cambridge Philos. Soc.} {\bf 53}, 629--641.\\
Cardy, J. (1992). Critical percolation in finite geometries. {\em  J. of Phy.  A } {\bf 25}  L201.\\
Grimmett, G.  (1999). {\em Percolation.} Springer-Verlag,  New York.\\
Kesten, H. (1982). {\em Percolation theory for mathematicians}, Birkhauser, Boston.\\
Langlands, R., Pouliot, P., and Saint-Aubin, Y. (1994). Conformal invariance in two-dimensional percolation. {\em Bull.  Amer. Math. Soc.} {\bf 30}  1--61.\\
Smirnov, S. (2001). Critical percolation in the plane: conformal invariance, Cardy's formula, scaling limits. {\em  C. R. Acad.  Sci. I-Math.} {\bf 333}  239--244.\\
Werner, W. (2008). Lectures on two-dimensional critical percolation, arXiv:0710.0856.\\
Yu Zhang\\
Department of Mathematics\\
University of Colorado\\
Colorado Springs, CO 80933\\
yzhang3@uccs.edu

\end{document}